\documentclass[12pt]{article}

\usepackage[margin=1in]{geometry}  
\usepackage{graphicx}              
\usepackage{amsmath}               
\usepackage{amsfonts}              
\usepackage{amsthm}                

\begin{document}

\nocite{*}

\title{Chow groups of moduli spaces of rank 2 vector bundles on curves with determinant of odd degree.}

\author{Evgeny Mayanskiy}

\maketitle

\begin{abstract}
  We compute Chow groups of moduli spaces of rank 2 vector bundles on curves with determinant of odd degree in terms of generators and relations.
\end{abstract}

\section{Theorem.}

{\bf Theorem. } {\it Let $X$ be a smooth projective curve of genus $g \geq 2$, $d=2w+1$, $w \geq 2g-2$, $\Lambda$ be a fixed line bundle on $X$ of degree $d$, and $N$ be the moduli space of semistable (= stable) rank $2$ vector bundles on $X$ of determinant $\Lambda$. Denote $n=d-2g+2$, $m=d+g-1$, $M_0=\mathbb P^{m-1}$ and $S^k X$ a symmetric power of $X$. All the Chow groups appearing below have integer coefficients.\\

Then for any $l$ there is a short exact sequence of abelian groups 
$$
0 \rightarrow \bigoplus_{i=0}^{n-2} A^{l+1+i}(N) \oplus \bigoplus_{k=1}^{w} \bigoplus_{r=0}^{k-2} A^{n-1+l-m+2k-r} (S^k X) \rightarrow
$$

$$
\rightarrow A^{n-1+l}(M_0) \oplus \bigoplus_{k=1}^{w} \bigoplus_{s=0}^{m-2k-2} A^{n-1+l-k-s} (S^k X) \rightarrow A^l (N) \rightarrow 0
$$
Moreover, there is a splitting
$$
\tau \colon A^l (N) \rightarrow A^{n-1+l}(M_0) \oplus \bigoplus_{k=1}^{w} \bigoplus_{s=0}^{m-2k-2} A^{n-1+l-k-s} (S^k X).
$$}

{\it Both the splitting $\tau$ and the morphisms in the exact sequence above are natural in the sense that they may be given explicitely, using standard intersection theory operations and some polynomials (whose coefficients can be determined by an inductive rule) in Chern classes of natural vector bundles defined on various $M_i$ and $S^k X$.}\\

{\bf Corollary. } {\it For any $l$ there is a finite resolution by abelian groups 
$$
0 \leftarrow A^l(N) \leftarrow \Omega_0 \leftarrow \Omega_1 \leftarrow \Omega_2 \leftarrow \ldots \leftarrow \Omega_t \leftarrow 0, 
$$
where each $\Omega_i$ is an (explicit) direct sum of Chow groups of $M_0$ and various $S^k X$, and the morphisms are natural in the above sense.}\\

\section{Proof.}

\subsection{Background and Notation.}

We deduce our theorem from Thaddeus' representation of (a projective bundle over) the moduli space $N$ as a sequence of blow-ups and blow-downs of a projective space $M_0=\mathbb P^{m-1}$. Let's recall the result due to Thaddeus, which we are going to use \cite{Thaddeus}.\\ 

There exists a projective bundle $M_w$ of rank $n$ over $N$, which can be put into the following sequence of blow-ups and blow-downs with smooth centers of smooth projective varieties:\\

\begin{figure}[ht]
\begin{center}
\includegraphics[scale=1]{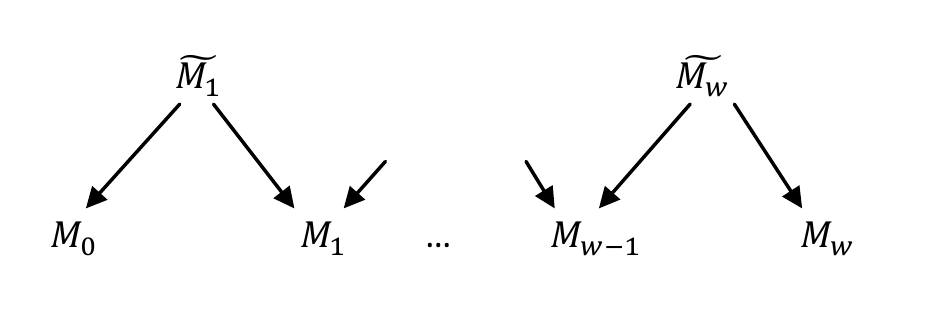}
\end{center}
\end{figure}

We use the following notations, when considering the $k-$th term in the Thaddeus' sequence:

\begin{figure}[ht]
\begin{center}
\includegraphics[scale=0.8]{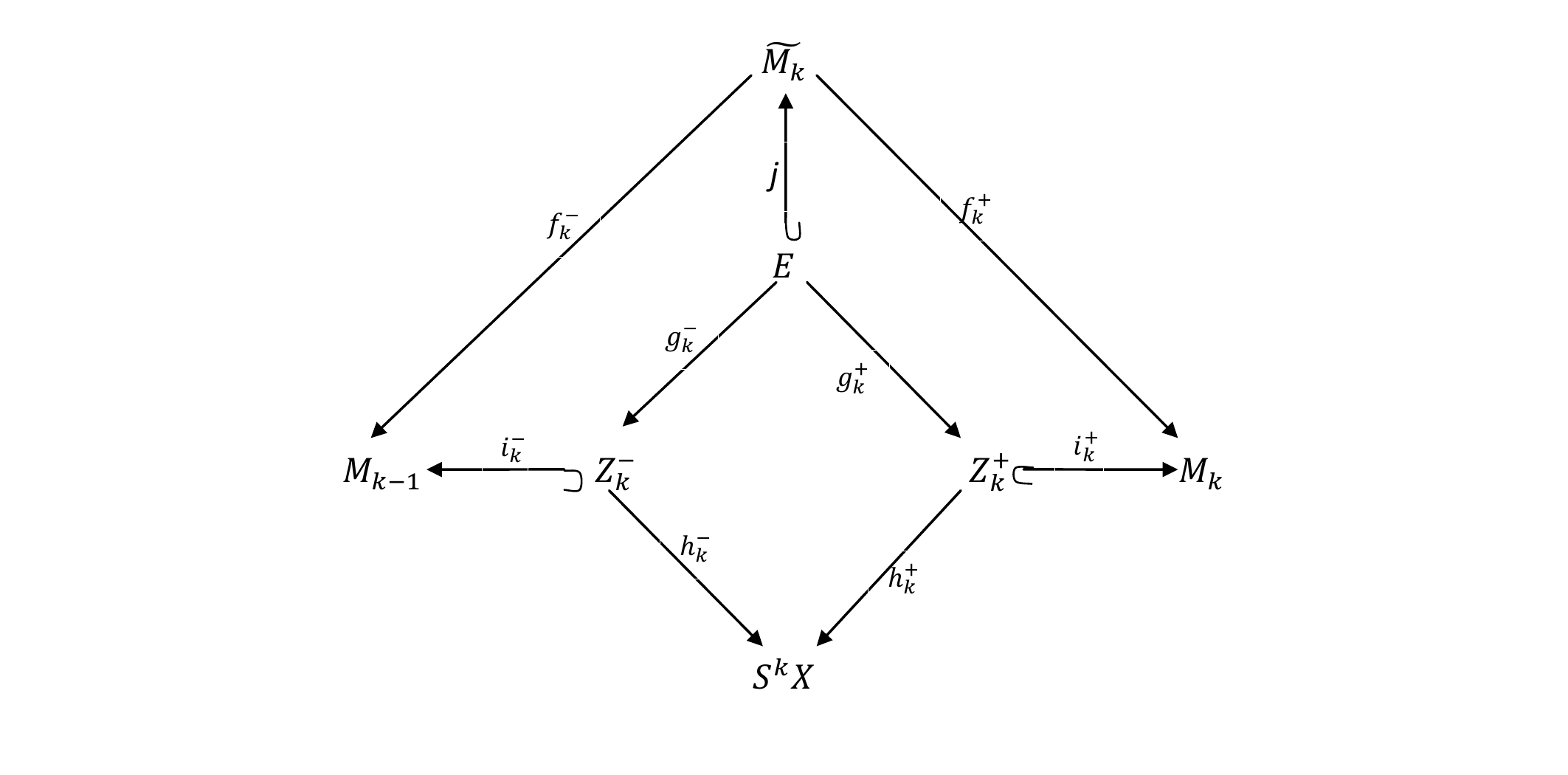}
\end{center}
\end{figure}

\begin{enumerate}
\item $W_{k}^{-}$ and $W_{k}^{+}$ are vector bundles of ranks $k$ and $m-2k$ on $S^kX$.

\item $Z_{k}^{-}=\mathbb PW_{k}^{-}$ is embedded by the map $i_{k}^{-}$ into $M_{k-1}$, $Z_{k}^{+}=\mathbb PW_{k}^{+}$ is embedded by the map $i_{k}^{+}$ into $M_{k}$.

\item $h_{k}^{\pm} \colon Z_{k}^{\pm} \rightarrow S^kX$ are the bundle maps, and ${\xi}_{k}^{\pm}$ be the first Chern classes of the corresponding line bundles $O(1)$.

\item $f_{k}^{-} \colon \tilde{M_{k}} \rightarrow M_{k-1}$ and $f_{k}^{+} \colon \tilde{M_{k}} \rightarrow M_{k}$ are blow-ups at $Z_{k}^{\pm}$ with exceptional divisor $E$, embedded by the map $j$ into $\tilde{M_k}$.
\item $g_{k}^{\pm} \colon E \rightarrow Z_{k}^{\pm}$ are the corresponding projections, which are projective bundles (since the central square is Cartesian) with the first Chern classes of the corresponding line bundles $O(1)$ denoted by ${\zeta}_{k}^{\pm}$ respectively.

\item Note that ${\zeta}_{k}^{+}=(g_{k}^{-})^{*} {\xi}_{k}^{-}$, ${\zeta}_{k}^{-}=(g_{k}^{+})^{*} {\xi}_{k}^{+}$, $i_{k}^{-}$ is a regular embedding of codimension $m-2k$, and $i_{k}^{+}$ is a regular embedding of codimension $k$.

\item $m=d+g-1$, $d=2w+1$, $w \geq 2g-2$.
\end{enumerate}

\subsection{Inductive computation of Chow groups of $M_i$'s}
We use the standard formulas describing Chow groups of a blow-up from \cite{Fulton}, section 6.7 twice. First, to express $A(\tilde{M_k})$ in terms of $A(M_{k-1})$ and some part of $A(E)$ (it will be a direct sum of these two terms). As an intermediate step, we use the formula for the projective bundle (twice), to compute $A(E)$ in terms of $A(S^kX)$. Finally, we use the formula for Chow groups of blow-ups again (for the blow-up on the right hand side) in order to express $A(M_k)$ in terms of $A(\tilde{M_k})$. Summarizing, this computes $A(M_k)$ in terms of $A(M_{k-1})$.\\

Let $f \colon \tilde Y \rightarrow Y$ be the blowup of $Y$ at $X$, which is embedded into $Y$ by $i \colon X \hookrightarrow Y$. Suppose $codim_Y(X)=d$. Let $g \colon \tilde X \rightarrow X$ be the exceptional divisor.\\

\begin{figure}[ht]
\begin{center}
\includegraphics[scale=0.8]{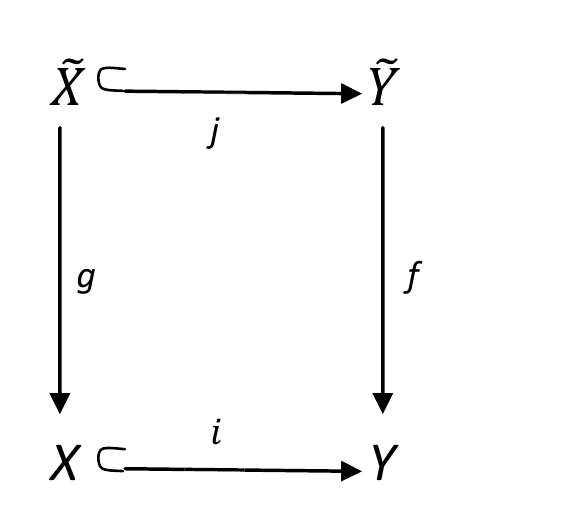}
\end{center}
\end{figure}

Let's denote by $\xi$ the first Chern class of $\mathcal E$, where $\tilde X={\mathbb P}_X (\mathcal E)$, and by $j \colon \tilde X \hookrightarrow \tilde Y$ the corresponding closed embedding. Then from Proposition 6.7 and Theorem 3.3 in \cite{Fulton}, we get an isomorphism:
$$
f^{*}(A^k(Y))\oplus j_{*}\left( \bigoplus_{s=0}^{d-2} \xi^s \cdot g^{*}(A^{k-s-1}(X)) \right) \rightarrow A^k(\tilde Y) 
$$

Consequently:
$$
0\rightarrow \bigoplus_{s=0}^{d-2} \xi^s \cdot g^{*}(A^{k-s-1}(X))\rightarrow A^k(\tilde Y)\rightarrow A^k(Y)\rightarrow 0
$$
where the first map send $\tilde x$ to $j_{*}\tilde x$, and the second map sends $\tilde y$ to $f_{*}(\tilde y)$.\\

Applying this to our situation, we get the split exact sequence:
$$
0\rightarrow \bigoplus_{s=0}^{k-2} (\zeta_k^{+})^s \cdot (g_k^{+})^{*} (A^{l-s-1}(Z_k^{+}))\xrightarrow{j_{*}} A^l(\tilde{M_k})=
$$

$$
=(f_k^{-})^{*}(A^l(M_{k-1}))\oplus j_{*} \left[ \bigoplus_{s=0}^{m-2k-2} ({\xi}_{k}^{-})^s\cdot (g_{k}^{-})^{*}(A^{l-s-1}(Z_{k}^{-})) \right]\xrightarrow{(f_{k}^{+})_{*}} A^l(M_k)\rightarrow 0
$$
with splitting:
$$
(f_k^{+})^{*} \colon A^l(M_k) \rightarrow A^l(\tilde{M_k})
$$
Recalling that $Z_k^{\pm}$ are projective bundles over $S^kX$, and applying Theorem 3.3 from \cite{Fulton}, we get:
$$
0 \rightarrow \bigoplus_{s=0}^{k-2} \bigoplus_{r=0}^{m-2k-1} ((\zeta)_k^{+})^s \cdot ((\zeta)_k^{-})^r \cdot (h_k\circ g_k)^{*}(A^{l-s-r-1}(S^kX))\xrightarrow{j_{*}}
$$

$$
\xrightarrow{j_{*}}(f_k^{-})^{*}(A^l(M_{k-1}))\oplus j_{*}\left[
\bigoplus_{s=0}^{m-2k-2} \bigoplus_{r=0}^{k-1} ((\zeta)_k^{-})^s \cdot ((\zeta)_k^{+})^r \cdot (h_k\circ g_k)^{*}(A^{l-s-r-1}(S^kX))
\right]\xrightarrow{(f_k^{+})_{*}} A^l(M_k)
\rightarrow 0
$$
with splitting given by $(f_k^{+})^{*} \colon A^l(M_k) \rightarrow \cdots$.\\

This sequence reduces to:
$$
0 \rightarrow \bigoplus_{r=0}^{k-2} ((\zeta)_k^{+})^r \cdot ((\zeta)_k^{-})^{m-2k-1} \cdot (h_k\circ g_k)^{*}(A^{l-m+2k-r}(S^kX))\xrightarrow{j_{*}}
$$

$$
\xrightarrow{j_{*}}(f_k^{-})^{*}(A^l(M_{k-1}))\oplus j_{*}\left[
\bigoplus_{s=0}^{m-2k-2} ((\zeta)_k^{-})^s \cdot ((\zeta)_k^{+})^{k-1} \cdot (h_k\circ g_k)^{*}(A^{l-k-s}(S^kX))
\right]\xrightarrow{(f_k^{+})_{*}} A^l(M_k)
\rightarrow 0
$$

This can be rewritten as:
$$
0 \longrightarrow \bigoplus_{r=0}^{k-2} (A^{l-m+2k-r}(S^kX)) \longrightarrow A^l(M_{k-1}) \oplus \bigoplus_{s=0}^{m-2k-2} A^{l-k-s} (S^k X) \longrightarrow A^l (M_k) \longrightarrow 0 \eqno{(*)}
$$
Maps here can be made explicit:
\begin{enumerate}
\item The map on the left sends $(v_r)_r$ to $\sum_{r=0}^{k-2} (i_{k}^{-})_{*} (({\xi}_{k}^{-})^r \cdot(h_{k}^{-})^{*}(v_r)) + (\sum_{r=0}^{k-2} C_{r}^{s} v_r)_s$.
\item The map on the right sends $m+(u_s)_s$ to $(f_{k}^{+})_{*} (f_{k}^{-})^{*}(m)+\sum_{s=0}^{m-2k-2} (i_{k}^{+})_{*} (({\xi}_{k}^{+})^s \cdot (h_{k}^{+})^{*}(u_s))$.
\item We used the following notation above: $C_{r}^{s} = \sum_{j=0}^{m-3k-s+r} (-1)^{j+r-k}\cdot {s+k-r+j \choose s+1} \cdot s_j(W_{k}^{-})\cdot c_{m-3k-s+r-j}(W_{k}^{+})$.
\end{enumerate}

This exact sequence has a splitting:
$$
A^l(M_k)\rightarrow A^l(M_{k-1})\oplus \bigoplus_{s=0}^{m-2k-2}A^{l-k-s}(S^kX) {\mbox , }\quad \quad n\mapsto (f_k^{-})_{*}(f_k^{+})^{*}(n)+(u_s)_s
$$
where $u_s$ can be given by a natural formula.\\

\subsection{Summing up.}
On the previous step we got a sequence of split exact sequences (with obvious notation):
$$
0 \longleftarrow A_k \longleftarrow A_{k-1} \oplus R_k \longleftarrow S_k \longleftarrow 0,
$$ 

$k=0,1,...,w$.\\

Because of the splitting, we can 'sum them up' to get the following split exact sequence:
$$
0 \longleftarrow A_w \longleftarrow A_{0} \oplus \bigoplus_{k=1}^{w} R_k \longleftarrow \bigoplus_{k=1}^{w} S_k \longleftarrow 0. \eqno{(**)}
$$ 

\subsection{Chow groups of $N$.}
Denote $n=d-2g+2$. We know that $M_w$ is a projective bundle of rank $n$ over $N$. Denoting by $\xi$ the first Chern class of the corresponding line bundle $O(1)$, we can write by the standard formula for Chow groups of projective bundles (Theorem 3.3 in \cite{Fulton}) the following split exact sequence:
$$
0 \longrightarrow \bigoplus_{i=0}^{n-2} (A^{l+n-1-i}(N)) \longrightarrow A^{n-1+l}(M_w) \longrightarrow A^l (N) \longrightarrow 0 \eqno{(***)}
$$
Now, summing this split exact sequence (as we did in the previous section) with sequence $(**)$, we get the exact sequence stated in the theorem. QED

\section{Acknowledgement.} I learned about this problem and about Thaddeus' sequence from Maxim Arap. 

\bibliographystyle{ams-plain}

\bibliography{ChowGroupsOfN.bib}

\end{document}